\documentclass[a4paper,11pt]{amsart}

\newtheorem{theorem}{Theorem}

\begin{document}

\title{Formal Adjoints and a Canonical Form for Linear Operators}

\renewcommand{\thefootnote}{}
\footnote{MGE is supported by the Australian Research Council. ARG is
supported by the New Zealand Institute for Mathematics and its
Applications and the Royal Society of New Zealand (Marsden Grant
02-UOA-108).}
\renewcommand{\subjclassname}{\textup{2000} Mathematics Subject Classification}
\renewcommand{\subjclassname}{\textup{2000} Mathematics Subject Classification}

\begin{abstract}
  We describe a canonical form for linear differential operators that are
formally self-adjoint or formally skew-adjoint.
\end{abstract}

\author{Michael G.\ Eastwood}
\address{Department of Pure Mathematics\\
University of Adelaide\\
South AUSTRALIA 5005
}\email{meastwoo@maths.adelaide.edu.au}

\author{A.\ Rod Gover}
\date{}
\address{Department of Mathematics\\
  The University of Auckland\\
  Private Bag 92019\\
  Auckland 1\\
  New Zealand}\email{gover@math.auckland.ac.nz}

\maketitle

Suppose $E$ and $F$ are smooth vector bundles on an oriented smooth
manifold~$M$. Let $vol$ denote the bundle of volume forms on~$M$.
The {\em formal adjoint\/} of a linear differential operator $L:E\to F$ is the
differential operator
$L^\ast:F^\ast\otimes vol\to E^\ast\otimes vol$
characterised by the equation
$$
\int_M\langle L^\ast\sigma,\tau\rangle =
  \int_M\langle\sigma,L\tau\rangle \quad\mbox{for
  }\sigma\in\Gamma(M,F^\ast\otimes vol)\mbox{ and }
  \tau\in\Gamma_\ast(M,E).
$$
Here all sections are assumed sufficiently smooth and $\Gamma_*$
indicates the space of compactly supported sections.  If
$F=E^\ast\otimes vol$, then $L^\ast:E\to F$ and there is a canonical
decomposition
$$L=L_++L_-=\textstyle{\frac{1}{2}}[L+L^\ast]+
            \textstyle{\frac{1}{2}}[L-L^\ast]$$
into self-adjoint and skew-adjoint parts.
Henceforth let us assume that $M$ is equipped with a preferred volume form and
a compatible torsion-free connection $\nabla$ (e.g.\ $M$ is Riemannian).

Suppose that $E$ is trivial. Then the preferred volume form
trivialises $F$, and so $L$ may be viewed as taking functions to functions
and written in terms of the given connection. A formula
for its adjoint is determined by integration by parts.  Suppose, for
example, that $L$ is second order. Then we may write
$$L=S^{ab}\nabla_a\nabla_b+ T^b\nabla_b + R,$$
where the tensor $S^{ab}$ is symmetric.
Adopting the convention that $\nabla_a$ acts on everything to its right, we can
re-express $L$ in the form
$$
L=\nabla_aS^{ab}\nabla_b+(\tilde{T}^b\nabla_b+\nabla_b\tilde{T}^b)+
\tilde{R}
$$
where the tensor $\tilde{T}$ and scalar field $\tilde{R}$ are given
by $\tilde{T}^b=\frac{1}{2}(T^b-(\nabla_aS^{ab}))$ and
$\tilde{R}=R-(\nabla_b\tilde{T}^b)$. This is congenial since clearly
$$L^\ast =\nabla_aS^{ab}\nabla_b-(\tilde{T}^b\nabla_b+\nabla_b\tilde{T}^b)+
      \tilde{R}.$$
In particular,
$$L_+=\nabla_aS^{ab}\nabla_b+\tilde{R}\quad\mbox{and}\quad
  L_-=\tilde{T}^b\nabla_b+\nabla_b\tilde{T}^b.$$ By an obvious
  inductive argument this generalises immediately to give the
  following result.
\begin{theorem}\label{one} A self-adjoint $k^{th}$ order linear
differential operator taking functions to functions on $M$ has even
order and may be canonically written in the form:
$$\sum_{i=0}^{k/2}\,\underbrace{\nabla_a \cdots \nabla_b}_i
  S_{(i)}^{a\cdots bc\cdots d}
          \underbrace{\nabla_c\cdots \nabla_d}_i,$$
for suitable symmetric tensors $S_{(i)}^{a\cdots d}$. A skew-adjoint
$k^{th}$ order linear differential operator taking functions to
functions on $M$ has odd order and may be
canonically written in the form:
$$\sum_{i=0}^{(k-1)/2}\underbrace{\nabla_a\cdots \nabla_b}_i
            (\nabla_c A_{(i)}^{a\cdots bcd\cdots e}
            +A_{(i)}^{a\cdots bcd \cdots e}\nabla_c)
        \underbrace{\nabla_d\cdots \nabla_e}_i,$$
for suitable symmetric tensors $A_{(i)}^{a\cdots e}$.
\end{theorem}

 Suppose now that $E$ is a vector bundle (possibly a tensor bundle)
with connection on $M$ and
write $\nabla$ to indicate the coupled tensor--vector bundle
connection. Let us use upper case Greek indices as abstract indices for
the bundle $E$ and its dual.  Then, for example, operators $L:E\to E^*$
may be written $ L_{\Psi\Phi}: E^\Phi\to E_\Psi $. Since the tensor
product $E^*\otimes E^*$ decomposes canonically into symmetric and
skew-symmetric parts it follows easily that the above generalises as
follows. We write $[\ell]$ to indicate the integer part of a real
number $\ell$.
\begin{theorem}[]
A self-adjoint (respectively skew-adjoint) $k^{th}$ order linear
differential operator $ L_{\Psi\Phi}: E^\Phi\to E_\Psi  $ on $M$
 may be canonically written in the form:
$$
\begin{array}{l}
\sum_{i=0}^{[k/2]}\,\underbrace{\nabla_a \cdots \nabla_b}_i
  S_{\Psi\Phi~(i)}^{a\cdots bc\cdots d}
          \underbrace{\nabla_c\cdots \nabla_d}_i+\\
\sum_{i=0}^{[(k-1)/2]}\underbrace{\nabla_a\cdots \nabla_b}_i
            (\nabla_c A_{\Psi\Phi~(i)}^{a\cdots bcd \cdots e}
            +A_{\Psi\Phi~(i)}^{a\cdots bcd\cdots e}\nabla_c)
        \underbrace{\nabla_d\cdots \nabla_e}_i
\end{array} ,
$$ where: the sections $S_{\Psi\Phi~(i)}^{a\cdots d}$ are
symmetric over the tensor indices and symmetric (respectively
skew-symmetric) over the pair $\Psi\Phi$ of vector bundle indices; the
sections $A_{\Psi\Phi~(i)}^{a\cdots e}$ are symmetric over
the tensor indices and skew-symmetric (respectively symmetric) over
the pair $\Psi\Phi$ of vector bundle indices.
\end{theorem}

Here is an example application. Suppose that $P$ is a self-adjoint
linear differential operator $P$ of order $k>0$ that takes functions
to functions on~$M$. Write $\tilde P$ for the modified operator
obtained by subtracting from $P$ the scalar part obtained by applying
$P$ to $f\equiv 1$. Then $ \tilde P$ is formally self-adjoint and
takes the form $f\mapsto \nabla_a(Q^{ab}(\nabla_bf))$ for a suitable
differential operator $Q: T^*M \to TM$. Since the construction of $\tilde P$
is canonical $\tilde P$ enjoys the same naturality and/or
invariance properties as does $P$.

A good example arises in conformal geometry.  Suppose $M$ is
an oriented conformal manifold of even dimension $2m$.  Let
$L$ denote any conformally invariant operator on functions that, with
respect to any Riemannian metric in the conformal class, takes the
form
\begin{equation}\label{Leq}
L=\Delta^{m}+\mbox{lower order terms}.
\end{equation}
 Since $L$ takes functions
to volume forms, so does its formal adjoint. The self-adjoint part $L_+$ of
$L$ is therefore also conformally invariant. As a conformal analogue
of $\Delta^{m}$ (in the sense of~\cite{ER}),
we may as well replace $L$ by~$L_+$. Then by Theorem
\ref{one} this, in turn, may be replaced by a self-adjoint conformally
invariant modification of the form
$f\mapsto\nabla_a(Q^{ab}(\nabla_bf))$ for a suitable $(n-2)^{nd}$
order differential operator $Q:\Lambda^1\to\Lambda^{n-1}$.
Thus, we obtain the following
result conjectured to us by Tom Branson.
\begin{theorem}\label{pline} Any conformally invariant linear differential
operator
 of the form {\rm(\ref{Leq})} admits a self-adjoint conformally invariant
modification of the form $f\mapsto\nabla_a(Q^{ab}(\nabla_bf))$ for a
suitable $(n-2)^{nd}$ linear order differential operator
$Q:\Lambda^1\to\Lambda^{n-1}$.
\end{theorem}
The motivation for his conjecture came from the case of the sphere
where the form of the operator may be verified directly. On the
sphere, the operator controls the embedding $L_{m}^2\hookrightarrow
e^L$ (Orlitz class) as a limiting case of the sharp Sobolev embeddings
$L_r^2\hookrightarrow L^{\frac{2m}{m-r}}$ for $r<n/2$ (equivalently,
comparing an $L^q$ norm with the complementary series
norm). See~\cite{b} for further discussion.  In \cite{gjms} Graham,
Jenne, Mason, and Sparling established the existence of operators of
the form (\ref{Leq}). More recently it has been shown, using
scattering theory~\cite{GrZ}, (and alternatively via the Fefferman-Graham ambient
metric \cite{BrGodeRham}) that the particular operators constructed in
\cite{gjms} are, in fact, already formally self-adjoint and have the
form given in Theorem~\ref{pline}.

\end{document}